\documentclass[12pt]{amsproc}
  \usepackage[margin=1in]{geometry}
\usepackage{setspace,fullpage}
\usepackage{graphicx}
\usepackage[nice]{nicefrac}
\usepackage{amssymb}
\usepackage{amsmath,amsthm,amscd,amssymb, mathrsfs}

\usepackage{latexsym}
\usepackage[colorlinks,citecolor=red,pagebackref,hypertexnames=false]{hyperref}
\usepackage{graphicx}
\usepackage{subfigure}
\usepackage{CJKutf8}

\numberwithin{equation}{section}

\theoremstyle{plain}
\newtheorem{theorem}{Theorem}[section]

\theoremstyle{definition}

\theoremstyle{remark}
\newtheorem{remark}[theorem]{Remark}

\newtheorem{case[theorem]}{Case}

\title[\parbox{14cm}{\centering{ An integral arising from dyadic average of Riesz transforms \hspace{1in}}} \quad]{An integral arising from dyadic average of Riesz transforms}

\author[]{Chih-Chieh Hung}
\address{Department of Mathematics\\ National Taiwan University\\
Taipei, 106 Taiwan}
\email{r07221003@ntu.edu.tw}

\author[C-Y. Shen]{Chun-Yen Shen}
\address{Department of Mathematics\\ National Taiwan University and National Center for Theoretical Sciences\\Taipei, 106 Taiwan}
\email{cyshen@math.ntu.edu.tw}

\thanks{supported by MOST through grant 108-2628-M-002-010-MY4}


\begin{document}
\begin{CJK}{UTF8}{bsmi}
\begin{abstract}
In the work of S. Petermichl, S. Treil and A. Volberg it was explicitly constructed that the Riesz transforms in any dimension $n \geq 2$ can be obtained as an average of dyadic Haar shifts provided that an integral is nonzero. It was shown in the paper that when $n=2$, the integral is indeed nonzero (negative) but for $n \geq 3$ the nonzero property remains unsolved. In this paper we show that the integral is nonzero (negative) for $n=3$. The novelty in our proof is the delicate decompositions of the integral for which we can either find their closed forms or prove an upper bound.  \end{abstract}
     
 \maketitle    
\tableofcontents

\section{Introduction} 
One of the important milestones that appears in the area of Harmonic analysis in the past decades is the appearance of dyadic Haar shifts. It not only connects the continuous singular integrals with dyadic operators but also enables people to resolve the longstanding $A_2$ conjecture concerning with the sharp weighted bound for Calder\'on-Zygmund singular integrals. More precisely, the breakthrough work of Petermichl \cite{Pt00} showed that the kernel of the Hilbert transform is actually an average of some certain dyadic operators:

\begin{equation}\label{hilbert}
\frac{c_0}{t-x}=\lim_{L\rightarrow\infty}\frac{1}{2\log L}\int^L_{\frac{1}{L}}\lim_{R\rightarrow\infty}\frac{1}{2R}\int^R_{-R} \Sigma_{I\in D^{\alpha,r}}h_I(t)(h_{I-}(x)-h_{I+}(x))d\alpha dr.
\end{equation}
Therefore the sharp weighted bound for Hilbert transform can be reduced to proving a uniform sharp weighted bound for above dyadic operators which are called dyadic Haar shifts. Such representation of dyadic average for Hilbert transform kernel later was generalized by Stefanie Petermichl, Sergei Treil, Alexander Volberg, \cite{Vo02} to a slightly wider class of kernels but still restricted on one dimensional singular integrals. Finally the full generality was made by Hyt\"onen in \cite{Hy3} who showed that any Calder\'on-Zygmund  operator is 
a simple variant of dyadic averages, and part of the work was built on a previous result obtained by Hyt\"onen, Perez, Treil and Volberg \cite{Hy2}.\\

Moreover as shown in the work of S. Petermichi that an explicit dyadic Haar shift can actually be given for Hilbert transform. As a result it may be also expected that some explicit dyadic Haar shifts can also be given for Riesz transforms. Indeed, it was shown in the work of Petermichl, Treil, Volberg \cite{Vo02} that  each component of the kernel of Riesz transforms can be explicitly represented by an average of dyadic shift upon showing the nonzero property of a constant through the process of averaging. Roughly speaking, it is in a sense equivalent to say that the constant $c_0$ that appears on the left hand side of above equation (\ref{hilbert}) is nonzero in order to make the result nondegenerate, and similarly the constant that appears through the process of averaging for Riesz transforms is given in the form below. Therefore a question that was risen in their work \cite{Vo02} is whether the following integral is zero or not (the detail definitions of some notations in this integral are given in next section):
\begin{align}\label{main integral}
\int_{S^{n-1}_+}<\omega,e_1>\xi_n(\omega)d\sigma(\omega).
\end{align}
They were able to show the integral is nonzero (in fact it is negative) when dimension $n=2$ but for dimension $n \geq 3$ the problem remains unsolved. \\

Therefore the purpose of this paper is to show the above integral when dimension $n=3$ is also negative. This was done via a careful and an efficient decomposition for the integral. For some terms in our decomposition we are able to show explicitly that their values are negative. For some other terms we are able to prove an upper bound. Combining all the estimates shows the integral is indeed negative. Now let us mention the difficulties of the integral for dimension $n=3$. First, the integrand functions in the integral are piecewise defined on some compact intervals, and the range of the integration is only half-sphere. Secondly, after we carefully analyse the integrand functions, one of the main difficulties then arises due to the mutual overlaps of their supports. More precisely, after using the sphere coordinates in the integral, the supports of the functions will create several difficulties since the behaviors of the points in these supports will now depend on the values of some complicated trigonometry functions. For these difficulties, it requires us to very carefully distinguish the range of the integrals in our decomposition. Finally for several integrals in our decomposition, we are able to show that their exact definite integrals can be computed. For the other integrals, we are not able to find their exact definite integrals, but we are able to find their upper bounds whose values can be explicitly estimated.

\section{Preliminary}
In this section we first introduce some notations that will be used frequently in this paper. Let $F_0$, $F$, $f_0$, and $f$ be defined as followings.
\begin{align*}F_0(x)&=\begin{cases}
    \frac{1}{2}-|x+\frac{1}{2}|         & -1\leq x<0 \\
    -( \frac{1}{2}-|x-\frac{1}{2}| )    & 0\leq x <1 \\
                                0    &\text{otherwise}
\end{cases};\\
F(x)&=F_0(x)-\frac{1}{2}F_0(x+1)-\frac{1}{2}F_0(x-1);\\
f_0(x)&=\begin{cases}
    1-|x|         & \text{ if } |x|\leq 1 \\
                                0       & \text{ otherwise}
\end{cases};\\
f(x)&=f_0(x)+\frac{1}{2}f_0(x+1)+\frac{1}{2}f_0(x-1).
\end{align*}
Note that $F$ is an odd function so that we may only describe $F$ on $x\geq 0$, i.e.
$$F(x)=\left\{\begin{aligned}
    \frac{-3}{4}+\frac{3}{2}|x-\frac{1}{2}|         & &if& &|x-\frac{1}{2}|\leq \frac{1}{2} \\
     \frac{1}{4}-\frac{1}{2}|x-\frac{3}{2}|         & &if& &|x-\frac{3}{2}|\leq \frac{1}{2} \\
                                0       & & & & x\geq 2
\end{aligned}\right.$$
and 
$f(x)=
    1-\frac{1}{2}|x|$ if $|x|\leq 2$ and $0$ otherwise.  See below for their graphs.

\begin{figure}[htbp]
\centering
\subfigure[$F_0$]{
\begin{minipage}[t]{0.4\linewidth}
\centering
\includegraphics[width=2in]{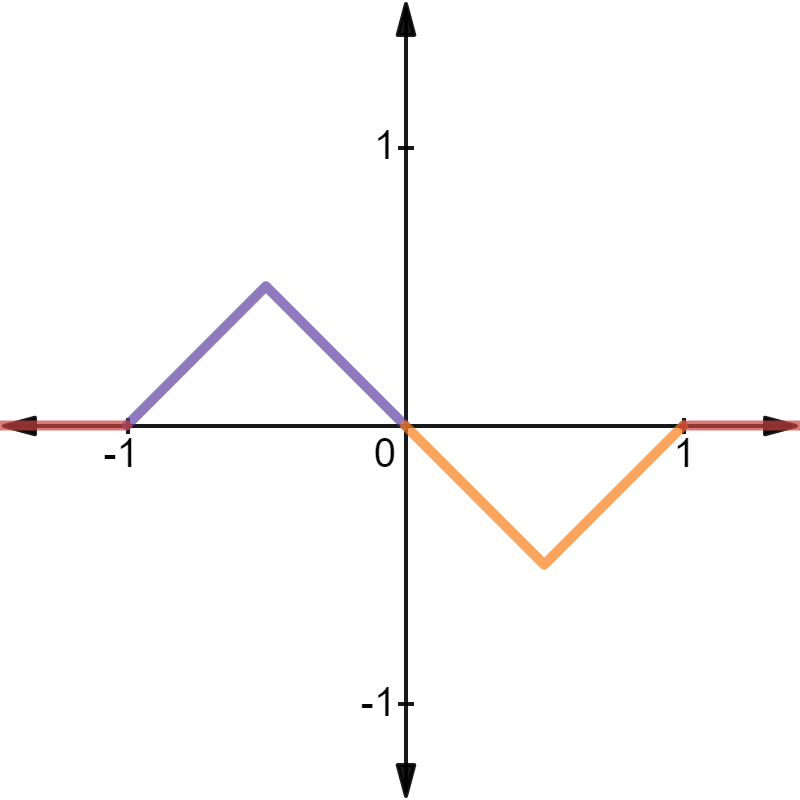}
\end{minipage}%
}%
\subfigure[$F$]{
\begin{minipage}[t]{0.4\linewidth}
\centering
\includegraphics[width=2in]{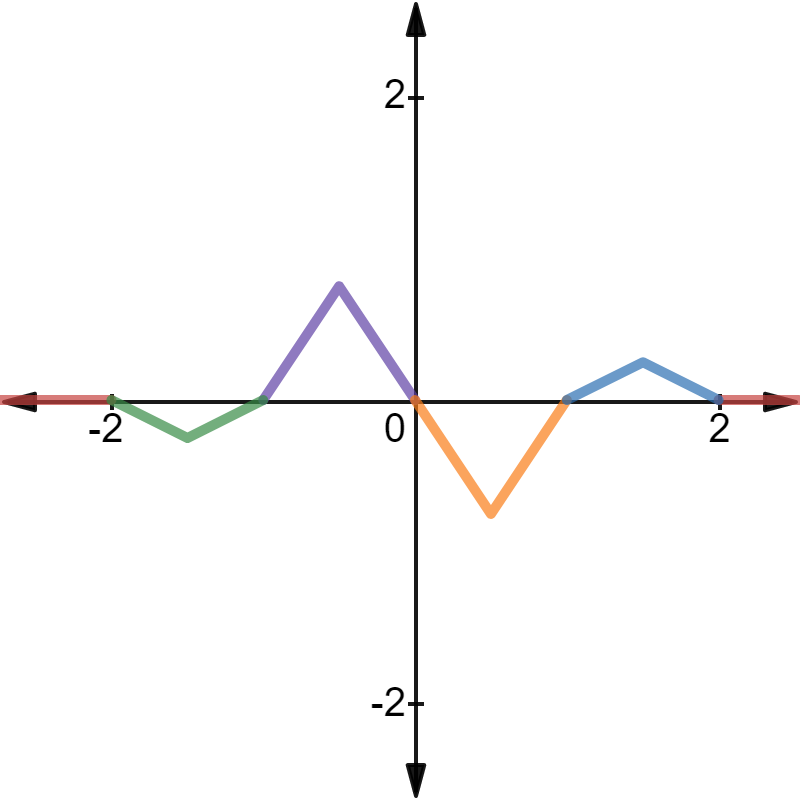}
\end{minipage}%
}%
\centering
\caption{ $F_0$ and $F$}
\end{figure}

\begin{figure}[htbp]
\centering
\subfigure[$f_0$]{
\begin{minipage}[t]{0.4\linewidth}
\centering
\includegraphics[width=2in]{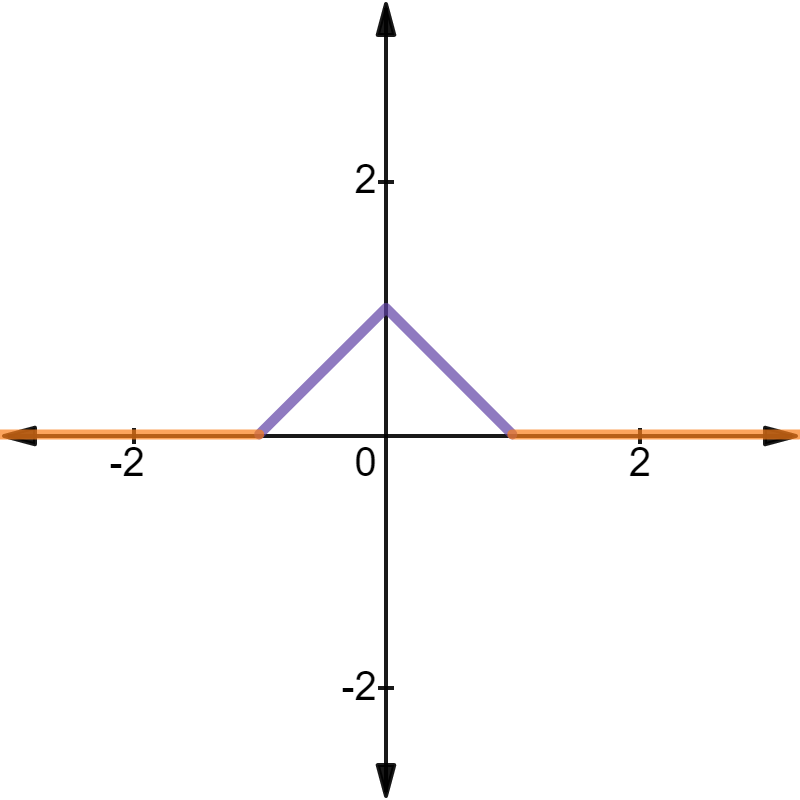}
\end{minipage}%
}%
\subfigure[$f$]{
\begin{minipage}[t]{0.4\linewidth}
\centering
\includegraphics[width=2in]{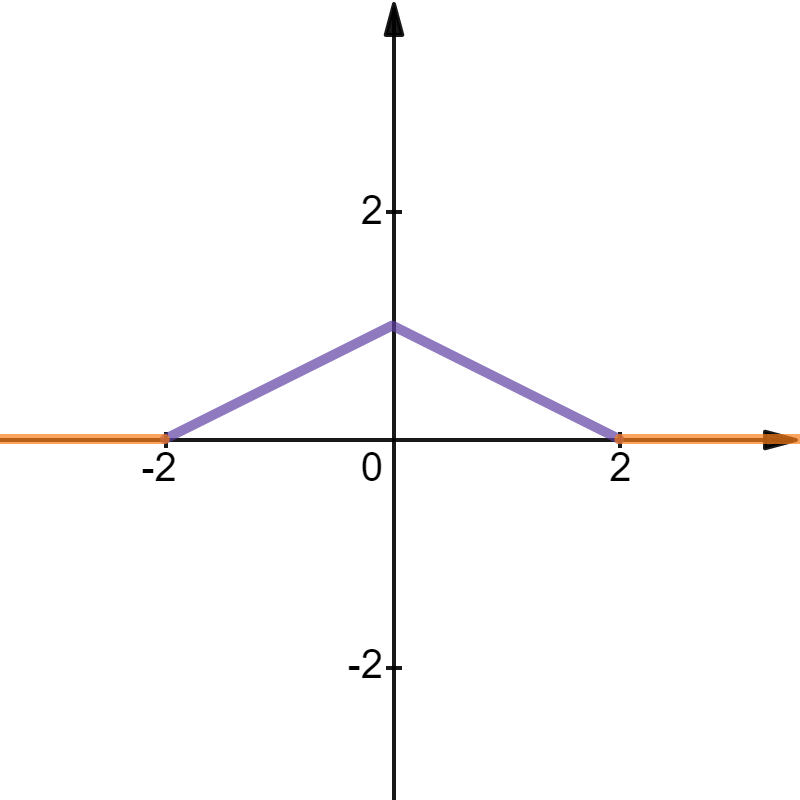}
\end{minipage}%
}%
\centering
\caption{ $f_0$ and $f$}
\end{figure}
\newpage
For all $n\geq 2$, and given $x=(x_1,...,x_n)\neq 0$, we define
\begin{align*}K_n(x)&=F(x_1)\times\prod_{i=2}^nf(x_i)\\
 \intertext{and}  
\xi_n(\frac{x}{|x|})&:=|x|^n\int_0^\infty\frac{1}{\rho^n}K_n(\frac{x}{\rho})\frac{d\rho}{\rho}.\end{align*}
Let $\rho=\frac{|x|}{t}$, the above formula $\xi_n(\frac{x}{|x|})$ now becomes
$$\xi_n(\omega):=\xi_n(\frac{x}{|x|})=\int_0^\infty t^{n-1}K_n(\frac{tx}{|x|})dt,$$
 where $\omega\in S^{n-1}$. Recall that the main result that we want to prove is to show the following integral is nonzero (negative),
 \begin{align}\label{main integral}
\int_{S^{n-1}_+}<\omega,e_1>\xi_n(\omega)d\sigma(\omega),
\end{align}
where $S^{n-1}_+=\{ \omega=(\omega_1, \cdots, \omega_n) \in S^{n-1} ; \omega_1 > 0 \}.$
 Thus putting it together, our goal is to show for $n=3$ the following integral is negative. 
 $$\int_{S^2_+}<\omega,e_1>\int_0^\infty t^2F(t\omega_1)f^{\omega_2}(t)f(t\omega_3)dtd\sigma(\omega).$$
 
 Before we proceed to give our decompositions for this integral. The diagram in next page gives a big picture how the integral is decomposed and how each term in our decomposition is estimated.
 
 \includegraphics[scale=0.8]{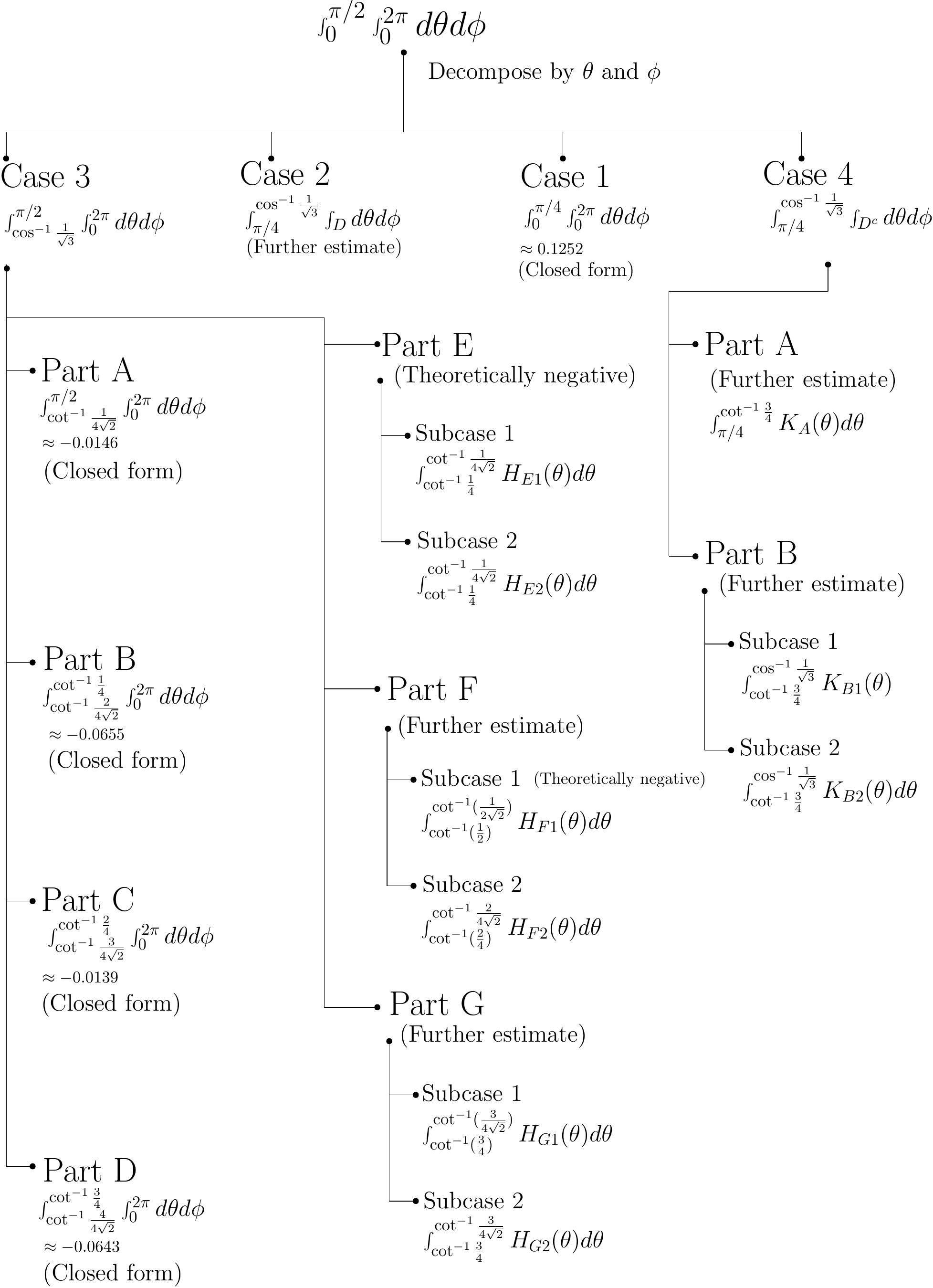}

\section{Decompositions}

In order to integrate with spherical measure appearing in (\ref{main integral}), we use spherical coordinate system to represent $\omega_1,\omega_2$, and $\omega_3$ i.e.
$$\omega_1=\cos\theta,$$
$$\omega_2=\sin\theta\cos\phi,$$
$$\omega_3=\sin\theta\sin\phi,$$
where $\theta\in [0,\pi/2]$ and $\phi\in [0,2\pi]$ because we only integrate on the half sphere $S^2_+$.
Putting in these new variables and using change of variables formula, the integral (\ref{main integral}) which we want to estimate becomes 
\begin{equation}\label{main integral3}
\int_0^{\frac{\pi}{2}}\int_0^{2\pi}\cos\theta\int_0^\infty t^2F(t\cos\theta)f(t\sin\theta\cos\phi)f(t\sin\theta\sin\phi)dt\sin\theta d\phi d\theta,
\end{equation}
where the factor $\sin\theta$ is due to Jacobian, and the $\cos\theta$ is from $<\omega,e_1>$. Since the integral range is $0\leq\theta\leq\frac{\pi}{2}$, we only need to consider $F^{\omega_1}(t)=F(t\cos\theta)$ with $\omega_1\geq 0$. In order to estimate $\xi_3(\omega)$, we break $F^{\omega_1}(t)$ into 4 linear mutually disjoint support functions with respect to $t$ that is $F^{\omega_1}(t):=F^{\omega_1}_{11}(t)+F^{\omega_1}_{12}(t)+F^{\omega_1}_{13}(t)+F^{\omega_1}_{14}(t)$, where 
\begin{equation}
F^{\omega_1}_{11}(t)=\left\{
\begin{array}{rcl}
-\frac{3}{2}t\omega_1& & {0\leq t\leq \frac{1}{2\omega_1}}\\
0 & & {\text{otherwise}}\\
\end{array} \right.
\end{equation}

\begin{equation}
F^{\omega_1}_{12}(t)=\left\{
\begin{array}{rcl}
\frac{3}{2}(t\omega_1-1) & & {\frac{1}{2\omega_1}\leq x\leq \frac{2}{2\omega_1}}\\
0 & & {\text{otherwise}}\\
\end{array} \right.
\end{equation}

\begin{equation}
F^{\omega_1}_{13}(t)=\left\{
\begin{array}{rcl}
\frac{1}{2}(t\omega_1-1) & & {\frac{2}{2\omega_1}\leq t\leq \frac{3}{2\omega_1}}\\
0 & & {\text{otherwise}}\\
\end{array} \right.
\end{equation}

\begin{equation}
F^{\omega_1}_{14}(t)=\left\{
\begin{array}{rcl}
-\frac{1}{2}(t\omega_1-2) & & {\frac{3}{2\omega_1}\leq t\leq \frac{4}{2\omega_1}}\\
0 & & {\text{otherwise}}\\
\end{array} \right.
\end{equation}

\begin{equation}
f^{\omega_2}(t)=\left\{
\begin{array}{rcl}
(1-\frac{|t\omega_2|}{2}) & & {0\leq |t|\leq \frac{2}{|\omega_2|}}\\
0 & & {\text{otherwise}}\\
\end{array} \right.
\end{equation}

\begin{equation}
f^{\omega_3}(t)=\left\{
\begin{array}{rcl}
(1-\frac{|t\omega_3|}{2}) & & {0\leq |t|\leq \frac{2}{|\omega_|}}\\
0 & & {\text{otherwise}}\\
\end{array} \right. .
\end{equation}

Let
$$S_i(a,\theta,\phi)=\int_0^a t^2 F^{\omega_1}_{1i}(t)f^{\omega_2}(t)f^{\omega_3}(t)dt=\int_0^a t^2 F^{\omega_1}_{1i}(t)f(t\sin\theta\cos\phi)f(t\cos\theta\sin\phi)dt,$$
for $i=1,2, 3, 4$. Then we have $\xi_3(\omega)=\sum_{i=1}^4 S_i(\infty,\theta,\phi)$. 
\begin{remark}
 For a fixed $\omega_1$, the support of $F^{\omega_1}(t)$ is $0\leq t\leq \frac{2}{|\omega_1|}$, and the support of $f^{\omega_2}$ and $f^{\omega_3}$ are $0\leq t\leq \frac{2}{|\omega_2|}$ and $0\leq t\leq \frac{2}{|\omega_3|}$ respectively. Hence the integral range for $S_i$ will be simultaneously determined by the supports of the functions of $F^{\omega_1}, f^{\omega_2}$, and $f^{\omega_3}$. This observation leads us to decompose the integral in terms of the supports of these functions. More precisely, we divide the integral into 4 cases depending on which function vanishes first. We now give details in the following sections.
\end{remark}

\subsection{Criteria of decomposition}
\text{}\\
As just mentioned before, we can reduce the integral range of $S_i$ to $\left[0,\max\{\frac{2}{|\omega_1|},\frac{2}{|\omega_2|},\frac{2}{|\omega_3|}\}\right]$. Therefore the integral range of $t$ now depends on the variables $\theta$ and $\phi$ since the variables $\omega_i$ depend on $\theta, \phi$ and this also explains why estimating the integral is complicated and difficult. Therefore we will have 4 cases that depend on which function $F^{\omega_1}$ $f^{\omega_2}$ or $f^{\omega_3}$ vanishes before the others. 

\begin{remark}
Throughout this paper, given two functions $f,g$ and assuming the supports of $f$ and $g$ are $[0,a]$ and $[0,b]$ respectively, then we say $f$ vanishes before $g$ or  $g$ vanishes after $f$ if $a\leq b$. 
\end{remark}
\subsection{$F^{\omega_1}(t)$ vanishes before the others.}  Assuming $F^{\omega_1}(t)$ vanishes before the others i.e. $|\omega_1|\geq |\omega_2|$ and $|\omega_1|\geq |\omega_3|$. Notice that $\omega_1^2+\omega_2^2+\omega_3^2=1$. If $|\omega_1|=|\cos\theta|\geq \frac{1}{\sqrt{2}}$, then for all $\phi$, $|\omega_1|$ is always the largest one. If $|\omega_1|\in [\frac{1}{\sqrt{3}},\frac{1}{\sqrt{2}}]$, then $ |\omega_1| \geq |\omega_2|$ for some $\phi$. If $|\omega_1|\leq \frac{1}{\sqrt{3}}$, then by pigeonhole principle one of $|\omega_2|, |\omega_3|$ must be larger than $\frac{1}{\sqrt{3}}$ which is larger than $|\omega_1|$. Thus $\theta$ must be in $[0,\cos^{-1}\frac{1}{\sqrt{3}}]$. There are two different situations we need to separately deal with. 
\begin{description}
\item[1] When $\theta\in [0,\pi/4],$ then we will have $\omega_1$ is larger than $|\omega_2|$ and $|\omega_3|$ for all $\phi$.
\item[2] When $\theta\in [\pi/4,\cos^{-1}\frac{1}{\sqrt{3}}],$ then  we will have $\omega_1$ is larger than $|\omega_2|$ and $|\omega_3|$ for some $\phi$.
\end{description}

\subsubsection{Case 1} $0\leq \theta< \pi/4$ and for all $0\leq \phi<2\pi$.\textbf{(Closed form)}

 Since $F^{\omega_1}(t)$ vanishes before the others and $\omega_1=\cos \theta$, $\omega_2= \sin\theta\cos\phi$ and $\omega_3=\sin\theta\sin\phi$. Thus we see that when $\theta\in [0,\pi/4],$ we will have
 $$|\omega_1|=|\cos \theta|\geq |\sin\theta\cos\phi|=|\omega_2|,$$ and $$|\omega_1|=|\cos\theta|\geq|\sin\theta\sin\phi|=|\omega_3|,$$
 for all $0\leq \phi<2\pi$.
 Now
 $$\xi_3(\omega)=\int^{\frac{2}{\omega_1}}_0t^2F(t\cos\theta)f(t\sin\theta\cos\phi)f(t\sin\theta\sin\phi)dt=\sum_{i=1}^4 S_i(\frac{i}{2\omega_1},\theta,\phi)-S_i(\frac{i-1}{2\omega_1},\theta,\phi)$$
let

\begin{align*}
    h_1(\theta,\phi)&:=S_1(\frac{1}{2\cos\theta},\theta,\phi),\\
    h_2(\theta,\phi)&:=S_2(\frac{1}{\cos\theta},\theta,\phi)-S_2(\frac{1}{2\cos\theta},\theta,\phi),\\
    h_3(\theta,\phi)&:=S_3(\frac{3}{2\cos\theta},\theta,\phi)-S_3(\frac{1}{\cos\theta},\theta,\phi),\\
    h_4(\theta,\phi)&:=S_4(\frac{2}{\cos\theta},\theta,\phi)-S_4(\frac{3}{2\cos\theta},\theta,\phi).\\
\end{align*}

Since $t^2F_{1i}(t\omega_1)f^{\omega_2}(t)f^{\omega_3}(t)$ is a polynomial of degree 5 on this integral range, and $h_1$, $h_2$, $h_3$, and $h_4$ all have closed forms, therefore we obtain 

\begin{align*}
    \sum_{i=1}^4h_i(\theta,\phi)=&-(720\cos\theta\sin\theta\cos\phi-680\cos^2\theta-629\cos\phi\sin\phi\\
    &+720\cos\theta\sin\theta\sin\phi+629\cos^2\theta\cos\phi\sin\phi)(3840\cos^5\theta)^{-1},
\end{align*}

\begin{align*}
    \int_0^{2\pi}\sin\theta\cos\theta\sum_{i=1}^4h_i(\theta,\phi)d\phi=\frac{4(\cos^2\theta(\frac{17\pi}{192}-\frac{629}{7680})-\frac{3\cos\theta\sin\theta}{8}+\frac{629}{7680})}{\cos^{5}\theta}.
\end{align*}

Therefore, we see that the integral below has a closed form so that we can estimate it accurately
$$\int_0^{\pi/4}\int_0^{2\pi}\sin\theta \cos\theta (h_1(\theta,\phi)+h_2(\theta,\phi)+h_3(\theta,\phi)+h_4(\theta,\phi)) d\phi d\theta \approx 0.1252.$$

\begin{remark}
This is the only term that has positive value.

\end{remark}

\subsubsection{Case 2} $\pi/4\leq \theta\leq \cos^{-1}{\frac{1}{\sqrt{3}}}$, and $\phi\in D$. \textbf{(further estimate)}\\
In this case $\pi/4\leq \theta\leq \cos^{-1}{\frac{1}{\sqrt{3}}}$ and since we want 
$\omega_1\geq |\omega_2|,$ and $\omega_1\geq |\omega_3|.$
Thus we have $\cot\theta\geq |\sin\phi|$ and $\cot\theta\geq |\cos\phi|$ in spherical coordinates. As a result, the range of $\phi$ will be restricted on $D$ which is 
\begin{align*}
  &\text{On the first quadrant:} &\cos^{-1}(\cot\theta)\leq &\phi\leq\sin^{-1}(\cot\theta) \\
  &\text{On the second quadrant:} &\frac{\pi}{2}+\cos^{-1}(\cot\theta)\leq &\phi\leq\frac{\pi}{2}+\sin^{-1}(\cot\theta) \\
  &\text{On the third quadrant:} &\pi+\cos^{-1}(\cot\theta)\leq &\phi\leq \pi +\sin^{-1}(\cot\theta)\\
  &\text{On the fourth quadrant:} &\frac{3 \pi}{2}+\cos^{-1}(\cot\theta)\leq &\phi\leq \frac{3 \pi}{2}+\sin^{-1}(\cot\theta).
\end{align*}
Observe that since $f$ is even so that for all $\phi$ we have
\begin{align*}
    F^{\omega_1}(t)f^{\omega_2}(t)f^{\omega_3}(t)&=F(t\cos\theta)f(t\sin\theta\cos\phi)f(t\sin\theta\sin\phi)\\
    &=F(t\cos\theta)f(t\sin\theta\cos(\phi+\pi/2))f(t\sin\theta\sin(\phi+\pi/2))\\
    &=F(t\cos\theta)f(t\sin\theta\sin\phi)f(t\sin\theta\cos\phi)\\
    &= F^{\omega_1}(t)f^{\omega_3}(t)f^{\omega_2}(t) .
\end{align*}
Thus 
\begin{align*}
    &\int_{\cos^{-1}(\cot\theta)}^{\sin^{-1}(\cot\theta)}(h_1(\theta,\phi)+h_2(\theta,\phi)+h_3(\theta,\phi)+h_4(\theta,\phi))d\phi\\
    =&\int_{\cos^{-1}(\cot\theta)+\frac{\pi}{2}}^{\sin^{-1}(\cot\theta)+\frac{\pi}{2}}(h_1(\theta,\phi)+h_2(\theta,\phi)+h_3(\theta,\phi)+h_4(\theta,\phi))d\phi\\
    =&\int_{\cos^{-1}(\cot\theta)+\frac{2\pi}{2}}^{\sin^{-1}(\cot\theta)+\frac{2\pi}{2}}(h_1(\theta,\phi)+h_2(\theta,\phi)+h_3(\theta,\phi)+h_4(\theta,\phi))d\phi\\
    =&\int_{\cos^{-1}(\cot\theta)+\frac{3\pi}{2}}^{\sin^{-1}(\cot\theta)+\frac{3\pi}{2}}(h_1(\theta,\phi)+h_2(\theta,\phi)+h_3(\theta,\phi)+h_4(\theta,\phi))d\phi.
\end{align*}
Hence the integral that we want to estimate is equal to
\begin{align}
    &\int_{\pi/4}^{\cos^{-1}\frac{1}{\sqrt{3}}}\cos\theta\int_{D}(h_1(\theta,\phi)+h_2(\theta,\phi)+h_3(\theta,\phi)+h_4(\theta,\phi))d\phi \sin\theta d\theta\label{case2}\\
    &=\int_{\pi/4}^{\cos^{-1}\frac{1}{\sqrt{3}}}4\cos\theta\int_{\cos^{-1}(\cot\theta)}^{\sin^{-1}(\cot\theta)}(h_1(\theta,\phi)+h_2(\theta,\phi)+h_3(\theta,\phi)+h_4(\theta,\phi))d\phi \sin\theta d\theta.\nonumber
\end{align}
Integrating with respect to $\phi$ is a closed form as above. However, when integrating with respect to $\theta$, we are unable to find its closed form. The reason is that the range is from $\cos^{-1}(\cot\theta)$ to $\sin^{-1}(\cot\theta)$, and after integrating the variable $\phi$ these upper and lower limits make the integrand in the variable $\theta$ extremely complicated. Therefore this case will be further estimated in the final section.

\subsection{$f^{\omega_2}$ or $f^{\omega_3}$ vanishes before the others.}

In this case it suffices to consider $f^{\omega_2}$ vanishes before the others. The reasons are the followings. First, we observe that  
\begin{align*}
     F^{\omega_1}(t)f^{\omega_2}(t)f^{\omega_3}(t)&=F(t\cos\theta)f(t\sin\theta\cos\phi)f(t\sin\theta\sin\phi)\\
    &=F(t\cos\theta)f(t\sin\theta\cos(\phi+\pi/2))f(t\sin\theta\sin(\phi+\pi/2)),
   \end{align*}
and notice that the if $f^{\omega_3}$ vanishes before the others, the integral range for $\phi$ in this case is only different from the integral range for $\phi$ by rotating $\frac{\pi}{2}$ in the case that $f^{\omega_2}$ vanishes before the others.   
Now we see that for $f^{\omega_2}$ vanishes before the others, we must have $ \phi \in [-\pi/4,\pi/4]$ and $ [-3 \pi/4,5 \pi/4]$.  Again since $f$ is even, it suffices to only consider the range $ [-\pi/4,\pi/4]$. Thus
\begin{align*}
    \int_0^{2\pi}\xi_3(\omega)d\phi&=\int_0^{2\pi}\int_0^\infty t^2 F^{\omega_1}(t)f^{\omega_2}(t)f^{\omega_3}(t)dt\\
    &=4\int_{-\pi/4}^{\pi/4}\int_0^\infty t^2 F^{\omega_1}(t)f^{\omega_2}(t)f^{\omega_3}(t)dtd\phi.
\end{align*}
We now determine the range of $\theta$. Assume $f^{\omega_2}$ vanishes before the others i.e. $|\omega_2|\geq|\omega_1|$ and $|\omega_2|\geq |\omega_3|$. Notice that $\omega_1^2+\omega_2^2+\omega_3^2=1$. If $|\omega_2|=|\sin\theta\cos\phi|\geq \frac{1}{\sqrt{2}}$, then for all $\phi$, $|\omega_2|$ is always the largest one. If $|\omega_2|\in [\frac{1}{\sqrt{3}},\frac{1}{\sqrt{2}}]$, then $|\omega_2|\geq |\omega_1|$ for some $\phi$. If $|\omega_2|\leq \frac{1}{\sqrt{3}}$, then by pigeonhole principle one of $|\omega_1|, |\omega_3|$ is larger than $\frac{1}{\sqrt{3}}$ which is larger than $|\omega_2|$. Thus $\theta$ must be in $[\pi/4,\pi/2]$.
 There are two different situations we need to separately deal with
\begin{description}
\item[1]When $\theta\in [\cos^{-1}\frac{1}{\sqrt{3}},\pi/2]$, $|\omega_2|$ is larger than $|\omega_1|$ and $|\omega_3|$ for all $\phi\in [-\pi/4,\pi/4]$ .
\item[2]When $\theta\in [\pi/4,\cos^{-1}\frac{1}{\sqrt{3}}]$, $|\omega_2|$ is larger than $|\omega_1|$ and $|\omega_3|$ for some $\phi\in [-\pi/4,\pi/4]$.
\end{description}

\subsubsection{Case 3} $\cos^{-1}(\frac{1}{\sqrt{3}})\leq \theta\leq \pi/2$.

We also break $F^{\omega_1}$ into 4 pieces as before. And the decompositions in this term are the most complicated one since we need to decide which of the 5 pieces $F^{\omega_1}_{11}$, $F^{\omega_1}_{11}+F^{\omega_1}_{12}$, $F^{\omega_1}_{11}+F^{\omega_1}_{12}+F^{\omega_1}_{13}$, $F^{\omega_1}_{11}+F^{\omega_1}_{12}+F^{\omega_1}_{13}+F^{\omega_1}_{14}$, and $f^{\omega_2}$ vanishes before the others according to $\theta$ and $\phi$.
More precisely, since $f^{\omega_2}(t)$ vanishes before the others and  $F(t\omega_1)=F^{\omega_1}_{11}+F^{\omega_1}_{12}+F^{\omega_1}_{13}+F^{\omega_1}_{14}$. Therefore there are 4 possibilities that $f^{\omega_2}$ vanishes before (1) $F^{\omega_1}_{11}$ (2) $F^{\omega_1}_{11}+F^{\omega_1}_{12}$, (3) $F^{\omega_1}_{11}+F^{\omega_1}_{12}+F^{\omega_1}_{13}$ and (4)  $F^{\omega_1}_{11}+F^{\omega_1}_{12}+F^{\omega_1}_{13}+F^{\omega_1}_{14}$. We now further decompose these 4 possibilities in terms of the range of $\phi$. \\

{\bf{1: These 4 possibilities hold for all $\phi \in [-\pi/4,\pi/4]$.}} \\

\begin{description}
\item[Part A]$f^{\omega_2}$ vanishes before $F^{\omega_1}_{11}$. It gives $$\xi_3(\omega)=\int_0^{\frac{2}{|\omega_2|}}t^2F^{\omega_1}(t)f^{\omega_2}(t)f^{\omega_3}(t)=\int_0^{\frac{2}{|\omega_2|}}t^2F_{11}^{\omega_1}(t)f^{\omega_2}(t)f^{\omega_3}(t).$$

\item[Part B]$f^{\omega_2}$ vanishes after $F^{\omega_1}_{11}$ and before $F^{\omega_1}_{11}+F^{\omega_1}_{12}$. It gives $$\xi_3(\omega)=\int_0^{\frac{2}{|\omega_2|}}t^2(F_{11}^{\omega_1}(t)+F^{\omega_1}_{12}(t))f^{\omega_2}(t)f^{\omega_3}(t).$$

\item[Part C]$f^{\omega_2}$ vanishes after $F^{\omega_1}_{11}+F^{\omega_1}_{12}$ and before $F^{\omega_1}_{11}+F^{\omega_1}_{12}+F^{\omega_1}_{13}$. It gives $$\xi_3(\omega)=\int_0^{\frac{2}{|\omega_2|}}t^2(F_{11}^{\omega_1}(t)+F^{\omega_1}_{12}(t)+F^{\omega_1}_{13}(t))f^{\omega_2}(t)f^{\omega_3}(t).$$

\item[Part D]$f^{\omega_2}$ vanishes after $F^{\omega_1}_{11}+F^{\omega_1}_{12}+F^{\omega_1}_{13}$ and before $F^{\omega_1}_{11}+F^{\omega_1}_{12}+F^{\omega_1}_{13}+F^{\omega_1}_{14}$. It gives $$\xi_3(\omega)=\int_0^{\frac{2}{|\omega_2|}}t^2(F_{11}^{\omega_1}(t)+F^{\omega_1}_{12}(t)+F^{\omega_1}_{13}(t)+F^{\omega_1}_{14}(t))f^{\omega_2}(t)f^{\omega_3}(t).$$
\end{description}
Now we define some notations to simplify our expressions. Let 
\begin{align*}
    g_1(\theta,\phi):&=S_1(\frac{2}{|\omega_2|},\theta,\phi)=S_1(\frac{2}{\sin\theta\cos\phi},\theta,\phi),\\
    g_2(\theta,\phi):&=S_2(\frac{2}{|\omega_2|},\theta,\phi)-S_2(\frac{1}{2\cos\theta},\theta,\phi)=S_2(\frac{2}{\sin\theta\cos\phi},\theta,\phi)-S_2(\frac{1}{2\cos\theta},\theta,\phi),\\
    g_3(\theta,\phi):&=S_3(\frac{2}{|\omega_2|},\theta,\phi)-S_3(\frac{2}{2\cos\theta},\theta,\phi)=S_3(\frac{2}{\sin\theta\cos\phi},\theta,\phi)-S_3(\frac{2}{2\cos\theta},\theta,\phi),\\
    g_4(\theta,\phi):&=S_4\left(\frac{2}{|\omega_2|},\theta,\phi\right)-S_4(\frac{3}{2\cos\theta},\theta,\phi)=S_4\left(\frac{2}{\sin\theta\cos\phi},\theta,\phi\right)-S_4(\frac{3}{2\cos\theta},\theta,\phi).
\end{align*}
Therefore
\begin{description}
\item[Part A] $\xi_3(\omega)=g_1(\theta,\phi).$
\item[Part B] $\xi_3(\omega)=h_1(\theta,\phi)+g_2(\theta,\phi).$
\item[Part C] $\xi_3(\omega)=h_1(\theta,\phi)+h_2(\theta,\phi)+g_3(\theta,\phi).$
\item[Part D] $\xi_3(\omega)=h_1(\theta,\phi)+h_2(\theta,\phi)+h_3(\theta,\phi)+g_4(\theta,\phi).$
\end{description}
Since we have further decomposed the integral into these 4 possibilities for all $\phi\in [-\pi/4,\pi/4]$, we need to determine the range of $\theta$.

\paragraph{Part A} \textbf{(Closed form)}

Since we now have $0\leq\frac{2}{|\omega_2|}\leq \frac{1}{2 |\omega_1|}$  or $0\leq\frac{2}{|\omega_3|}\leq \frac{1}{2|\omega_1|}$ for all $\phi\in [-\pi/4,\pi/4]$. 
 In order to satisfy the condition $\max\{\frac{2}{|\omega_2|},\frac{2}{|\omega_3|}\}\leq\frac{1}{2|\omega_1|}$, thus we have $\max\{\frac{2}{|\omega_2|},\frac{2}{|\omega_3|}\}\leq\frac{2\sqrt{2}}{\sin\theta}\leq\frac{1}{2\cos\theta}=\frac{1}{2\omega_1}$, which in turn gives that
$$\pi/2\geq\theta\geq \cot^{-1}\frac{1}{4\sqrt{2}} .$$
 Therefore we get the following integral:
$$H_A(\theta):=4\sin\theta\cos\theta\int_{-\pi/4}^{\pi/4}g_1(\theta,\phi)d\phi=\frac{8\sin^2\theta-1}{\sin^3\theta}.$$
Moreover the integral of $H_A(\theta)$ with respect to $\theta$ is a closed form. Finally we can explicitly compute the vaule 
$$\int^{\pi/2}_{\cot^{-1}\frac{1}{4\sqrt{2}}}H_A(\theta)d\theta\approx -0.0146.$$

\paragraph{Part B} \textbf{(Closed form)}

Since we have $\frac{1}{2\omega_1}\leq\frac{2}{|\omega_2|}\leq\frac{2}{2\omega_1}$ for all $\phi$. In other words we have
$\frac{2\sqrt{2}}{\sin\theta}\leq \frac{2}{2\cos\theta}$, and $\frac{1}{2\cos\theta}\leq \frac{2}{\sin\theta}$. Hence we get range of $\theta$
$$\cot^{-1}\frac{1}{4}\geq\theta\geq \cot^{-1}\frac{2}{4\sqrt{2}}.$$
Therefore
\begin{align*}
    H_B(\theta):&=4\sin\theta\cos\theta\int_{-\pi/4}^{\pi/4}h_1(\theta,\phi)+g_2(\theta,\phi)d\phi\\
    &=-(24\cos\theta\sin^5\theta-\sin^6\theta-1024\cos^6\theta+2048\cos^5\theta\sin\theta-40\pi\sin^4\theta+40\pi\sin^6\theta\\
    &+5120\log(\sqrt{2}+1)\cos^5\theta\sin\theta+1024\sqrt{2}\cos^5\theta\sin\theta)/(1280\cos^4\theta\sin^3\theta).
\end{align*}
Similarly, the integral of $H_B(\theta)$ is a closed form and finally we have
$$\int_{\cot^{-1}\frac{2}{4\sqrt{2}}}^{\cot^{-1}\frac{1}{4}}H_B(\theta)d\theta\approx -0.0655.$$

\paragraph{Part C} \textbf{(Closed form)}

Since we have:
$\frac{2}{2\omega_1}\leq\frac{2}{|\omega_2|}\leq\frac{3}{2\omega_1}$ for all $\phi\in [-\pi/4,\pi/4]$. As above, we will have that the range of $\theta$ is
$$\cot^{-1}\frac{2}{4}\geq\theta\geq \cot^{-1}\frac{3}{4\sqrt{2}}.$$
Therefore
\begin{align*}
    H_C(\theta):&=4\sin\theta\cos\theta\int_{\pi/4}^{\pi/4}h_1(\theta,\phi)+h_2(\theta,\phi)+g_3(\theta,\phi)d\phi\\
    &=-[61\sin\theta-696\cos\theta+2088\cos^3\theta+\cos^5\theta(5120\log(\sqrt{2}+1)+1024\sqrt{2}-40)\\
    &-\cos^7\theta(5120\log(\sqrt{2}+1)+1024\sqrt{2}+1352)+(520\pi-183)(\sin\theta-\sin^3\theta)\\
    &-\cos^4\theta\sin\theta(1040\pi-183)+\cos^6\theta\sin\theta(520\pi-10301)/(3840\cos^4\theta\sin^4\theta)],
\end{align*}
and finally we have
$$\int_{\cot^{-1}\frac{3}{4\sqrt{2}}}^{\cot^{-1}\frac{2}{4}}H_C(\theta)d\theta\approx -0.0139.$$

\paragraph{Part D} \textbf{(Closed form)}

Since we have:
$\frac{3}{2\omega_1}\leq\frac{2}{|\omega_2|}\leq\frac{4}{2\omega_1}$ for all $\phi\in [-\pi/4,\pi/4]$. As above, we will have that the range of $\theta$ is
$$\cot^{-1}\frac{3}{4}\geq\theta\geq \cot^{-1}\frac{4}{4\sqrt{2}}.$$
Therefore
\begin{align*}
    H_D(\theta):&=4\sin\theta\cos\theta\int_{\pi/4}^{\pi/4}h_1(\theta,\phi)+h_2(\theta,\phi)+h_3(\theta,\phi)+g_4(\theta,\phi)d\phi\\
    &=-(395\sin\theta-3264\cos\theta+9792\cos^3\theta + \cos^7\theta(5120\log(\sqrt{2} + 1) + 1024\sqrt{2} + 5312)\\
    &-\cos^5\theta(5120\log(\sqrt{2} + 1)+ 1024\sqrt{2}+ 11840) + (1880\pi - 1185)(\sin\theta - \sin^3\theta)\\
    &- \cos^4\theta\sin\theta(3760\pi - 1185) + cos^6\theta\sin\theta(1880\pi + 4725))/(1920\cos^4\theta\sin^4\theta)),
\end{align*}
and finally we have
$$\int_{\cot^{-1}\frac{4}{4\sqrt{2}}}^{\cot^{-1}\frac{3}{4}}H_D(\theta)d\theta\approx -0.0634$$

However there are still some ranges of $\theta$ that we have not dealt with in parts $A, B, C, D$ above and the ranges are $$[\cot^{-1}\frac{1}{4},\cot^{-1}\frac{1}{4\sqrt{2}}],\,\,\,\,\,[\cot^{-1}\frac{2}{4},\cot^{-1}\frac{2}{4\sqrt{2}}],\text{ and}\,\,\,\,\,\left[\cot^{-1}\frac{3}{4},\cot^{-1}\frac{3}{4\sqrt{2}}\right].$$  
For these  3 ranges, we need to further estimate the integrals.

\begin{remark}
Those 3 parts are more complicated than above 4 parts since in each cases, $\sup\{\textit{supp} f^{\omega_2}\}=\frac{2}{|\omega_2|}$ is not contained in one of $[0,\frac{1}{2\omega_1}]$, $[\frac{1}{2\omega_1},\frac{2}{2\omega_1}]$, $[\frac{2}{2\omega_1},\frac{3}{2\omega_1}]$, and $[\frac{3}{2\omega_1},\frac{4}{2\omega_1}]$ for all $\phi\in[-\pi/4,\pi/4]$. In fact it depend on $\phi$. Therefore those three parts are not part A,B,C and D. Fortunately, it is contained in one of $[0,\frac{2}{2\omega_1}]$, $[\frac{1}{2\omega_1},\frac{3}{2\omega_1}]$, and $[\frac{2}{2\omega_1},\frac{4}{2\omega_1}]$. That is why in each of following cases, we need to split the integrals into two subcases according to $\phi$.
\end{remark}

\paragraph{Part E} $\theta\in [\cot^{-1}\frac{1}{4},\cot^{-1}\frac{1}{4\sqrt{2}}]$.\textbf{ (negative value)}

When $\theta\in [\cot^{-1}\frac{1}{4},\cot^{-1}\frac{1}{4\sqrt{2}}]$ the range $[- \pi/4, \pi/4]$ of $\phi$ will be split into two cases. One of the ranges will give $\xi_3(\omega)=g_1(\theta,\phi)$, and the other range will give $\xi_3(\phi)=g_2(\theta,\phi)$. More precisely, when  $\theta\in [\cot^{-1}\frac{1}{4},\cot^{-1}\frac{1}{4\sqrt{2}}]$, the variable $\phi$ will have two possibilities. One possibility is that $f^{\omega_2}$ vanishes before $F^{\omega_1}_{11}$ i.e. we will have $0 \leq\frac{2}{|\omega_2|}\leq\frac{1}{|2\omega_1|}$. The other possibility is that $f^{\omega_2}$ vanishes before $F^{\omega_1}_{12}$ and after $F^{\omega_1}_{11}$, which gives $\frac{1}{2\omega_1}\leq \frac{2}{|\omega_2|}\leq\frac{2}{2\omega_1}$. As a result, we split the range of $\phi$ according to which above possibility occurs.

\paragraph{{\bf{1}}}$0 \leq\frac{2}{|\omega_2|}\leq\frac{1}{2\omega_1}$ ($\xi_3(\omega)=g_1(\theta,\phi)$).

Since $0 \leq\frac{2}{|\omega_2|}\leq\frac{1}{2\omega_1}$ which is 
$\frac{2}{\sin\theta\cos\phi}\leq\frac{1}{2\cos\theta}$. Thus
$$-\cos^{-1}4\cot\theta\leq\phi\leq \cos^{-1}4\cot\theta.$$ Therefore we can estimate the integral below
$$H_{E1}(\theta):=4\sin\theta\cos\theta\int_{-\cos^{-1}(4\cot\theta)}^{\cos^{-1}(4\cot\theta)}g_1(\theta,\phi)d\phi,$$
which is a closed form in variable $\phi$. However after plugging in the upper and lower limits, we are unable to show whether the integral $\int H_{E1}(\theta) d \theta $ has a closed form. But it is easy to show that its value is negative. Notice that 
$g_1(\theta,\phi)=\int_0^{4/2\sin\theta\cos\phi}t^2F^{\omega_1}_{11}(t)f^{\omega_2}(t)f^{\omega_3}(t)dt$, and for all $t$, $t^2$, $f^{\omega_2}(t)$, and $f^{\omega_2}(t)$ are positive but $F^{\omega_1}_{11}(t)$ is negative so that $g_1(\theta,\phi)\leq 0$ for all $\theta$ and $\phi$. Hence
$$\int_{\cot^{-1}\frac{1}{4}}^{\cot^{-1}\frac{1}{4\sqrt{2}}}H_{E1}(\theta)d\theta\leq 0.$$
\paragraph{\bf{2}} $\frac{1}{2\omega_1}\leq \frac{2}{|\omega_2|}\leq \frac{2}{2\omega_1}$ ($\xi_3(\omega)=h_1+g_2$).

The integral range of $\phi$ is just the complement of the range in case {\bf{1}} above. Hence we have
$$H_{E2}(\theta):=4\sin\theta\cos\theta\left[\int_{\cos^{-1}(4\cot\theta)}^{\pi/4}h_1(\theta,\phi)+g_2(\theta,\phi)d\phi+\int^{-\cos^{-1}(4\cot\theta)}_{-\pi/4}h_1(\theta,\phi)+g_2(\theta,\phi)d\phi\right].$$
Again we can explicitly compute $H_{E2}(\theta)$ because the above integrals are closed forms in variable $\phi$. However after plugging in the upper and lower limits, the  integral $\int H_{E2}(\theta) d\theta$ is difficult to see if it has a closed form. But it is easy to show that its value is negative. Notice that
$h_1(\theta,\phi)+g_2(\theta,\phi)=\int_0^{4/2\sin\theta\cos\phi}t^2\left[F^{\omega_1}_{11}(t)+F^{\omega_1}_{12}(t)\right]f^{\omega_2}(t)f^{\omega_3}(t)dt$, and for all $t$, $t^2$, $f^{\omega_2}(t)$, and $f^{\omega_2}(t)$ are positive but $F^{\omega_1}_{11}(t)$, and $F^{\omega_1}_{12}(t)$ are negative so that $h_1(\theta,\phi)+g_2(\theta,\phi)\leq 0$ for all $\theta$ and $\phi$. Hence
$$\int_{\cot^{-1}\frac{1}{4}}^{\cot^{-1}\frac{1}{4\sqrt{2}}}H_{E2}(\theta)d\theta\leq 0.$$
\bigskip
\paragraph{Part F} $\theta\in [\cot^{-1}\frac{2}{4},\cot^{-1}\frac{2}{4\sqrt{2}}] $\textbf{(further estimate)}

As in part E above, there will be two cases when $\theta\in [\cot^{-1}\frac{2}{4},\cot^{-1}\frac{2}{4\sqrt{2}}].$ One case is that  $\frac{1}{2\omega_1}\leq\frac{2}{|\omega_2|}\leq\frac{2}{2\omega_1}$, and the other case is $\frac{2}{2\omega_1}\leq\frac{2}{|\omega_2|}\leq\frac{3}{2\omega_1}$.

\paragraph{\bf{1}}$\frac{1}{2\omega_1}\leq\frac{2}{|\omega_2|}\leq\frac{2}{2\omega_1}$ ($\xi_3(\omega)=h_1+g_2$).

Since $\frac{1}{2\omega_1}\leq\frac{2}{|\omega_2|}\leq\frac{2}{2\omega_1}$, it gives that $0\leq\phi\leq\cos^{-1}(2\cot\theta)$. Hence
$$H_{F1}(\theta):=4\sin\theta\cos\theta\int^{\cos^{-1}(2\cot\theta)}_{-\cos^{-1}(2\cot\theta)}h_1(\theta,\phi)+g_2(\theta,\phi)d\phi.$$
$H_{F1}(\theta)$ can be explicitly computed, but integral  $ \int H_{F1}(\theta) d\theta$ is difficult to show if it has a closed form. Therefore this case will be further estimated in the final section.

\paragraph{\bf{2}}$\frac{2}{2\omega_1}\leq\frac{2}{|\omega_2|}\leq\frac{3}{2\omega_1}$ ($\xi_3(\omega)=h_1+h_2+g_3$)

Since $\frac{2}{2\omega_1}\leq\frac{2}{|\omega_2|}\leq\frac{3}{2\omega_1}$. It gives us that 
\begin{align*}
    H_{F2}(\theta):=4\sin\theta\cos\theta(&\int^{\pi/4}_{\cos^{-1}(2\cot\theta)} h_1(\theta,\phi)+h_2(\theta,\phi)+g_3(\theta,\phi)d\phi\\
    &+\int_{-     \pi/4}^{-\cos^{-1}(2\cot\theta)} h_1(\theta,\phi)+h_2(\theta,\phi)+g_3(\theta,\phi)d\phi).
\end{align*}

Again $H_{F2}(\theta)$ can be explicitly computed, but integral  $\int H_{F2}(\theta) d\theta$ is difficult to show if it has a closed form. Therefore this case will be further estimated in the final section.
\bigskip
\paragraph{Part G} $\theta\in[\cot^{-1}\frac{3}{4},\cot^{-1}\frac{3}{4\sqrt{2}}]$ \textbf{(further estimate)}

\paragraph{\bf{1}}$\frac{2}{2\omega_1}\leq\frac{2}{|\omega_2|}\leq\frac{3}{2\omega_1}$ ($\xi_3(\omega)=h_1+h_2+g_3$).

Since $\frac{2}{2\omega_1}\leq\frac{2}{|\omega_2|}\leq\frac{3}{2\omega_1}$, it gives us that
$$H_{G1}(\theta):=4\sin\theta\cos\theta\int_{-\cos^{-1}(\frac{3}{4}\cot\theta)}^{\cos^{-1}(\frac{3}{4}\cot\theta)}h_1(\theta,\phi)+h_2(\theta,\phi)+g_3(\theta,\phi)d\phi.$$
$H_{G1}(\theta)$ can be explicitly computed, but integral $\int H_{G1}(\theta)$  is difficult to show if it has a closed form. Therefore this case will be further estimated in the final section.

\paragraph{\bf{2}}$\frac{3}{2\omega_1}\leq\frac{2}{|\omega_2|}\leq\frac{4}{2\omega_1}$ ($\xi_3(\omega)=h_1+h_2+h_3+g_4$)
Since $\frac{3}{2\omega_1}\leq\frac{2}{|\omega_2|}\leq\frac{4}{2\omega_1}$, it gives us that
\begin{align*}
    H_{G2}(\theta):=4\sin\theta\cos\theta(\int^{\pi/4}_{\cos^{-1}(\frac{3}{4}\cot\theta)} h_1(\theta,\phi)+h_2(\theta,\phi)+h_3(\theta,\phi)+g_4(\theta,\phi)d\phi\\
    +\int_{-\pi/4}^{-\cos^{-1}(\frac{3}{4}\cot\theta)} h_1(\theta,\phi)+h_2(\theta,\phi)+h_3(\theta,\phi)+g_4(\theta,\phi)d\phi)
\end{align*}
$$H_{G2}(\theta):=\int^{\pi/4}_{\cos^{-1}(\frac{3}{4}\cot\theta)} h_1(\theta,\phi)+h_2(\theta,\phi)+h_3(\theta,\phi)+g_4(\theta,\phi)d\phi.$$
$H_{G2}(\theta)$ can be explicitly computed, but integral $\int H_{G2}(\theta)$  is difficult to show if it has a closed form. Therefore this case will be further estimated in the final section.\\

Finally the remaining case is below.
\bigskip
\subsubsection{Case 4}: $\frac{\pi}{4}\leq \theta\leq\cos^{-1}(\frac{1}{\sqrt{3}})$ \textbf{(further estimate)}

For this case $\frac{\pi}{4}\leq \theta\leq\cos^{-1}(\frac{1}{\sqrt{3}})$, the range of $\phi$ is actually the complement of the range in case 2 above. In other words, the integral range of $\phi$ is $D^c$ (see page 8 for the definition of $D$). Just as what we observe in case 2 the integral that we want to estimate can be reduced to  

\begin{align}
    &\int_{\frac{\pi}{4}}^{\cos^{-1}\frac{1}{\sqrt{3}}}\cos\theta\int_0^{2\pi}t^2F^{\omega_1}(t)f^{\omega_2}(t)f^{\omega_3}(t)dtd\phi\sin\theta d\theta\nonumber\\
    &=\int_{\frac{\pi}{4}}^{\cos^{-1}\frac{1}{\sqrt{3}}}\cos\theta\int_{D^c}t^2F^{\omega_1}(t)f^{\omega_2}(t)f^{\omega_3}(t)dtd\phi\sin\theta d\theta\label{case4}\\
    &=\int_{\frac{\pi}{4}}^{\cos^{-1}\frac{1}{\sqrt{3}}}4\cos\theta\int_{-\cos^{-1}\cot\theta}^{\cos^{-1}\cot\theta}t^2F^{\omega_1}(t)f^{\omega_2}(t)f^{\omega_3}(t)dtd\phi\sin\theta d\theta.\nonumber
\end{align}

 Now here is the key observations that since $\phi\in [-\cos^{-1}\cot\theta,\cos^{-1}\cot\theta]$ so that $$\omega_2=\sin\theta\cos\phi\in [\cos\theta,\sin\theta],$$
 Hence we have 
$$\frac{2}{|\sin\theta|}\leq\frac{2}{|sin\theta\cos\phi|}\leq\frac{2}{|\cos\theta|}=\frac{2}{|\omega_1|}.$$
Therefore we will only have $f^{\omega_2}$ vanishes before $F^{\omega_1}_{11}+F^{\omega_1}_{12}+F^{\omega_1}_{13}+F^{\omega_1}_{14}$ and after $F^{\omega_1}_{11}+F^{\omega_1}_{12}+F^{\omega_1}_{13}$ for all $\phi\in [-\cos^{-1}\cot\theta,\cos^{-1}\cot\theta]$. This gives us that
$\frac{2}{|\sin\theta|}\geq \frac{3}{2|\omega_1|}$ and hence we have
$$\frac{\pi}{4}\leq\theta\leq \cot^{-1}\frac{3}{4}.$$
Therefore the other part is $$\cot^{-1}\frac{3}{4}\leq\theta\leq \cos^{-1}\frac{1}{\sqrt{3}}.$$

\paragraph{Part A'}: $\frac{\pi}{4}\leq\theta\leq \cot^{-1}\frac{3}{4}$\\
In this case we have
$$K_{A}(\theta):=4\cos\theta\sin\theta\int^{\cos^{-1}\cot\theta}_{-\cos^{-1}\cot\theta} h_1(\theta,\phi)+h_2(\theta,\phi)+h_3(\theta,\phi)+g_4(\theta,\phi)d\phi,$$
which shows we need to estimate the integral below for this case. 
$$\int_{\pi/4}^{\cot^{-1}\frac{3}{4}}K_{A}(\theta)d\theta.$$
Still, we are not able to find its closed form and hence this case will be further estimated in the final section.
\bigskip
\paragraph{Part B'}: $\cot^{-1}\frac{3}{4}\leq\theta\leq \cos^{-1}\frac{1}{\sqrt{3}}$

When $\cot^{-1}\frac{3}{4}\leq\theta\leq \cos^{-1}\frac{1}{\sqrt{3}}$, there will be two cases\\
(1) $f^{\omega_2}(t)$ vanishes after $F^{\omega_1}_{11}+F^{\omega_1}_{12}+F^{\omega_1}_{13}$, which gives $\frac{2}{\omega_2}\geq \frac{3}{2\omega_1}$;\\
(2) $f^{\omega_2}(t)$ vanishes before $F^{\omega_1}_{11}+F^{\omega_1}_{12}+F^{\omega_1}_{13}$ after $F^{\omega_1}_{11}+F^{\omega_1}_{12}$, which gives $\frac{2}{|\omega_2|}\leq \frac{3}{2\omega_1}$.

\paragraph{\bf{1}} $\frac{2}{|\omega_2|}\geq \frac{3}{2\omega_1}$ ($\xi_3(\omega)=h_1+h_2+h_3+g_4$).

From this condition, we get the integral range of $\phi$ is $$\cos^{-1}\frac{4}{3}\cot\theta\leq\phi\leq \cos^{-1}\cot\theta.$$
Therefore we have
\begin{align*}
    K_{B1}(\theta):=4\cos\theta\sin\theta(&\int^{\cos^{-1}\cot\theta}_{\cos^{-1}\frac{4}{3}\cot\theta}h_1(\theta,\phi)+h_2(\theta,\phi)+h_3(\theta,\phi)+g_4(\theta,\phi)d\phi\\
    &+\int_{-\cos^{-1}\cot\theta}^{-\cos^{-1}\frac{4}{3}\cot\theta}h_1(\theta,\phi)+h_2(\theta,\phi)+h_3(\theta,\phi)+g_4(\theta,\phi)d\phi.
\end{align*}
$K_{B1}$ can be explicitly computed, but integral $\int K_{B1}(\theta)$ is difficult to see if it has a closed form. Therefore this case will be further estimated in the final section.

\paragraph{\bf{2}} $\frac{2}{2\omega_1}\leq\frac{2}{|\omega_2|}\leq\frac{3}{2\omega_1}$ ($\xi_3(\omega)=h_1+h_2+g_3$).

From this condition, we get the integral range of $\phi$ is
$$0\leq\phi\leq \cos^{-1}\frac{4}{3}\cot\theta.$$

$$K_{B2}(\theta):=4\cos\theta\sin\theta\int^{\cos^{-1}\frac{4}{3}\cot\theta}_{-\cos^{-1}\frac{4}{3}\cot\theta}h_1(\theta,\phi)+h_2(\theta,\phi)+g_3(\theta,\phi)d\phi.$$
$K_{B2}(\theta)$ can be explicitly computed, but integral $\int K_{B2}$ is difficult to see if it has a closed form. Therefore this case will be further estimated in the final section.
 
\section{Some Remarks}
(1) We use the program of Matlab \cite{Mat} to find the closed forms for some of above integrals. It can be also directly checked that all the indefinite integrals are correct.  (2) The variable-precision floating-point arithmetic (VPA) that we use in the program of Matlab is 32 digits, thus the precision of the values for closed forms is accurate up to $10^{-32}$ error which would not effect our final value.

\section{Upper bounds for all the further estimate cases}
Recall that for all the closed forms above, their values add up to be negative. In addition, part E is proved to be negative. Therefore our goal in the section is to give upper bounds for all the further estimate cases above and show that the values of the upper bounds are all negative which after all shows the integral \ref{main integral} is negative.

Note that $\xi_3(\omega)=\int^\infty_0 t^2F(t\omega_1)f(t\omega_2)f(t\omega_3)dt$ and $F(t\omega_1)=F^{\omega_1}_{11}(t)+F^{\omega_1}_{12}(t)+F^{\omega_1}_{13}(t)+F^{\omega_1}_{14}(t)$, also recall that for all $t,\omega$, $F^{\omega_1}_{11}(t), F^{\omega_1}_{12}(t)$ are negative, and $F^{\omega_1}_{13}(t), F^{\omega_1}_{14}(t), f(t\omega_2), f(t\omega_3)$ are positive. In previous sections, we have shown that for those further estimate cases, it is difficult to see if they have closed forms. Therefore the ideas for these remaining cases are to combine the integrals and split the combined integrals into positive and negative integrals. Finally, we are able to find upper bounds for these negative and positive integrals and show that these upper bounds have closed forms. 
\bigskip
\subsection{Case 2 and case 4}
Recall that in case 2, the integral (\ref{case2}) is
$$\int_{\pi/4}^{\cos^{-1}\frac{1}{\sqrt{3}}}\cos\theta\sin\theta\int_{\phi\in D}\xi_3(\omega)d\phi  d\theta.$$
In case 4, the integral (\ref{case4}) is
$$\int_{\pi/4}^{\cos^{-1}\frac{1}{\sqrt{3}}}\cos\theta\sin\theta\int_{\phi\in D^c}\xi_3(\omega)d\phi  d\theta.$$
We now separate $\xi_3(\omega)$ into negative part and positive part i.e. 
\begin{align*}
    \textit{(Negative)}&&\int_{\pi/4}^{\cos^{-1}\frac{1}{\sqrt{3}}}\cos\theta\sin\theta\int_0^{2\pi}\int_0^{2\sqrt{3}}t^2
    \left[F^{\omega_1}_{11}(t)+F^{\omega_1}_{12}(t)\right]f^{\omega_2}(t)f^{\omega_3}(t)dtd\phi  d\theta, \\
    \textit{(Positive)}&&\int_{\pi/4}^{\cos^{-1}\frac{1}{\sqrt{3}}}\cos\theta\sin\theta\int_0^{2\pi}\int_0^{2\sqrt{3}}t^2
    \left[F^{\omega_1}_{13}(t)+F^{\omega_1}_{14}(t)\right]f^{\omega_2}(t)f^{\omega_3}(t)dtd\phi  d\theta,
\end{align*}
where negative part indicates that $t^2
    \left[F^{\omega_1}_{11}(t)+F^{\omega_1}_{12}(t)\right]f^{\omega_2}(t)f^{\omega_3}(t)\leq 0$ for all $t$ and positive part indicates $t^2
    \left[F^{\omega_1}_{13}(t)+F^{\omega_1}_{14}(t)\right]f^{\omega_2}(t)f^{\omega_3}(t)\geq 0$ for all $t.$                                                          
    
\subsubsection{\textbf{Negative part}}

Since in case 2 and case 4 the three functions $F, f^{\omega_2}, f^{\omega_3}$ all vanish after $F^{\omega_1}_{12}$ so that the negative part becomes
$$\int_{\pi/4}^{\cos^{-1}\frac{1}{\sqrt{3}}}\cos\theta\sin\theta\int_0^{2\pi}h_1(\theta,\phi)+h_2(\theta,\phi)d\phi d\theta.$$
Also this integral has closed form
\begin{align*}
    &\int\cos\theta\sin\theta\int_0^{2\pi}h_1(\theta,\phi)+h_2(\theta,\phi)d\phi d\theta\\
    &=\left[\frac{7\pi}{32} -\frac{9}{64}\tan\frac{\theta}{2} -\frac{280\pi-31}{640}\tan^2\frac{\theta}{2}+\frac{7\pi}{32}\tan^4\frac{\theta}{2} +\frac{9}{64}\tan^5\frac{\theta}{2}-\frac{31}{1920}\right]/(\tan^2\frac{\theta}{2}- 1)^3\\
    &-\frac{9\tanh^{-1}(\tan\frac{\theta}{2})}{64}. 
\end{align*}
Plugging in the exact integral range, we thus obtain
 $$\int_{\pi/4}^{\cos^{-1}\frac{1}{\sqrt{3}}}\cos\theta\sin\theta\int_0^{2\pi}h_1(\theta,\phi)+h_2(\theta,\phi)d\phi d\theta\approx -0.0607.$$
Thus the negative part of case 2 +case 4 is about -0.0607.

\subsubsection{\textbf{Positive part}}
For positive part, we are unable to get the exact value. Instead, we will find an upper bound for the positive part and show that the upper bound has a closed form. Recall that in case 2, and part A' and part B'1 of case 4, the function $t^2F(t\omega_1)f(t\omega_2)f(t\omega_3)$ vanishes after $F^{\omega_1}_{13}$ and they are
\begin{align*}
    &\textit{Case 2}
    &\int_{\pi/4}^{\cos^{-1}\frac{1}{\sqrt{3}}}\cos\theta\sin\theta\int_D h_3(\theta,\phi)+h_4(\theta,\phi)d\phi d\theta,\\
    &\textit{Case 4,A'}
    &\int_{\pi/4}^{\cot^{-1}\frac{3}{4}}\cos\theta\sin\theta\int_{D^c} h_3(\theta,\phi)+g_4(\theta,\phi)d\phi d\theta,\\
    &\textit{Case 4,B'1}
    &\int_{\cot^{-1}\frac{3}{4}}^{\cos^{-1}\frac{1}{\sqrt{3}}}\cos\theta\sin\theta\int_{R} h_3(\theta,\phi)+g_4(\theta,\phi)d\phi d\theta,
\end{align*}
where $R$ is the integral range of $\phi$ in B'1 of case 4.

Notice that in part B'2 of case 4, the function $t^2F(t\omega_1)f(t\omega_2)f(t\omega_3)$ vanishes before $F^{\omega_1}_{14}$ and the integral is
\begin{align*}
    &\textit{Case 4,B'2} &\int_{\cot^{-1}\frac{3}{4}}^{\cos^{-1}\frac{1}{\sqrt{3}}}\cos\theta\sin\theta\int_{R^c} g_3(\theta,\phi)d\phi d\theta.
\end{align*}

To obtain an upper bound, we introduce a function which is a linear extension of $F^{\omega_1}_{13}$.

\begin{equation*}
    \Tilde{F}^{\omega_1}_{13}(t)=\frac{1}{2}(t\omega_1-1).
\end{equation*}
\begin{center}
    \includegraphics[scale=0.35]{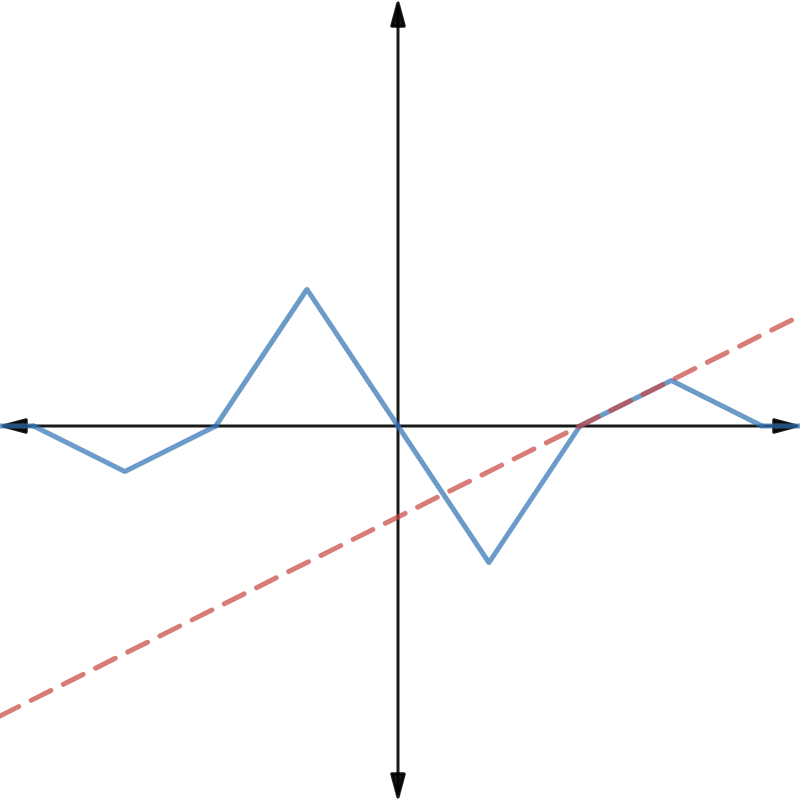}
\end{center}

Where the dotted line is $\Tilde{F}^{\omega_1}_{13}(t)=\frac{1}{2}(t\omega_1-1).$

Now for case 2
\begin{align*}
    h_3(\theta,\phi)+h_4(\theta,\phi)=&\int_{2/2\cos\theta}^{3/2\cos\theta}t^2F^{\omega_1}_{13}(t)f^{\omega_2}(t)f^{\omega_3}(t)dt+\int_{3/2\cos\theta}^{4/2\cos\theta}t^2F^{\omega_1}_{14}(t)f^{\omega_2}(t)f^{\omega_3}(t)dt\\
    \leq&\int_{2/2\cos\theta}^{3/2\cos\theta}t^2\Tilde{F}^{\omega_1}_{13}(t)f^{\omega_2}(t)f^{\omega_3}(t)dt+\int_{3/2\cos\theta}^{4/2\cos\theta}t^2\Tilde{F}^{\omega_1}_{13}(t)f^{\omega_2}(t)f^{\omega_3}(t)dt\\
    \leq&\int_{2/2\cos\theta}^{4/2\sin\theta\cos\phi}t^2\Tilde{F}^{\omega_1}_{13}(t)f^{\omega_2}(t)f^{\omega_3}(t)dt.
\end{align*}
The last inequality comes from that fact that $F$ vanishes before $f^{\omega_2}$ in case 2, i.e. 
$$\frac{4}{2\cos\theta}\leq\frac{4}{2\sin\theta\cos\phi}.$$
For case 4A and 4B'1
\begin{align*}
    h_3(\theta,\phi)+g_4(\theta,\phi)=&\int_{2/2\cos\theta}^{3/2\cos\theta}t^2F^{\omega_1}_{13}(t)f^{\omega_2}(t)f^{\omega_3}(t)dt+\int_{3/2\cos\theta}^{4/2\sin\theta\cos\phi}t^2F^{\omega_1}_{14}(t)f^{\omega_2}(t)f^{\omega_3}(t)dt\\
    \leq&\int_{2/2\cos\theta}^{3/2\cos\theta}t^2\Tilde{F}^{\omega_1}_{13}(t)f^{\omega_2}(t)f^{\omega_3}(t)dt+\int_{3/2\cos\theta}^{4/2\sin\theta\cos\phi}t^2\Tilde{F}^{\omega_1}_{13}(t)f^{\omega_2}(t)f^{\omega_3}(t)dt\\
    =&\int_{2/2\cos\theta}^{4/2\sin\theta\cos\phi}t^2\Tilde{F}^{\omega_1}_{13}(t)f^{\omega_2}(t)f^{\omega_3}(t)dt.
\end{align*}
For the case 4B'2, we note that
\begin{align*}
    g_3(\theta,\phi)&=\int_{2/2\cos\theta}^{3/2\sin\theta\cos\phi}t^2 F^{\omega_1}_{13}(t)f(t\omega_2)f(t\omega_3)dt\\
    &\leq \int_{2/2\cos\theta}^{4/2\sin\theta\cos\phi}t^2 \Tilde{F}^{\omega_1}_{13}(t)f(t\omega_2)f(t\omega_3)dt.
\end{align*}

Therefore combining all the inequalities above to get an upper bound which is
\begin{align*}
    &\int\cos\theta\sin\theta\int_0^{2\pi}\int_{2/2\cos\theta}^{4/2\sin\theta\cos\phi}t^2\Tilde{F}^{\omega_1}_{13}(t)f^{\omega_2}(t)f^{\omega_3}(t)dtd\phi d\theta\\
    =&(\tan^6\frac{\theta}{2}\frac{\pi+2}{6}- \tan^4\frac{\theta}{2}\frac{10\pi+29}{30}+ \tan^2\frac{\theta}{2}\frac{15\pi+89}{90} + \tan\frac{\theta}{2}(\frac{2\log(\sqrt{2}+1)}{3} + \frac{2\sqrt{2}+4}{15})\\
    &+ \tan^5\frac{\theta}{2}(2\log(\sqrt{2}+1) + \frac{2\sqrt{2}+4}{5}) - \tan^3\frac{\theta}{2}(2\log(\sqrt{2}+1) + \frac{4\sqrt{2}+9}{10}) \\
    &- \tan^7\frac{\theta}{2}(\frac{2\log(\sqrt{2}+1)}{3} + \frac{2\sqrt{2}}{15}+\frac{1}{6}) - \frac{1}{3})/(\tan^2\frac{\theta}{2} - 3\tan^4\frac{\theta}{2} + 3\tan^6\frac{\theta}{2} - \tan^8\frac{\theta}{2})\\
    &- (4\log(\tan\frac{\theta}{2}))/3 - \tanh^{-1}(\frac{6391}{9\tan\frac{\theta}{2} + 240)} - \frac{80}{3})/10 + \tan\frac{\theta}{2}(\frac{2\log(\sqrt{2}+1)}{3} + \frac{2\sqrt{2}+4}{15})\\
    &+ \frac{\tan\frac{\theta}{2}^2}{3}.
\end{align*}
Plugging in the exact integral range, we thus obtain 
$$\int_{\pi/4}^{\cos^{-1}\frac{1}{\sqrt{3}}}\cos\theta\sin\theta\int_0^{2\pi}\int_{2/2\cos\theta}^{4/2\sin\theta\cos\phi}t^2\Tilde{F}^{\omega_1}_{13}(t)f^{\omega_2}(t)f^{\omega_3}(t)dtd\phi d\theta\approx 0.08718.$$
Thus the positive part of case 2 +case 4 is bounded by 0.08717.
\subsection{Part F in case 3} Part $F$ is
\begin{align*}
    &\int^{\cot^{-1}\frac{2}{4\sqrt{2}}}_{\cot^{-1}\frac{2}{4}}4\sin\theta\cos\theta\int^{\cos^{-1}(2\cot\theta)}_{-\cos^{-1}(2\cot\theta)}h_1(\theta,\phi)+g_2(\theta,\phi)d\phi d\theta\\
    &\int^{\cot^{-1}\frac{2}{4\sqrt{2}}}_{\cot^{-1}\frac{2}{4}}4\sin\theta\cos\theta(\int^{\pi/4}_{\cos^{-1}(2\cot\theta)} h_1(\theta,\phi)+h_2(\theta,\phi)+g_3(\theta,\phi)d\phi d\theta\\
     &+\int^{\cot^{-1}\frac{2}{4\sqrt{2}}}_{\cot^{-1}\frac{2}{4}}\int_{-     \pi/4}^{-\cos^{-1}(2\cot\theta)} h_1(\theta,\phi)+h_2(\theta,\phi)+g_3(\theta,\phi)d\phi)d \theta,
\end{align*}
where $h_1,h_2,g_2$ represent the negative part and $g_3$ represents the positive part.
\subsubsection{\textbf{Negative part of F}} We first notice that the negative part of $F$ is bounded by
\begin{align*}
   &\int4\cos\theta\sin\theta\int_{-\pi/4}^{\pi/4}h_1(\theta,\phi)d\phi d\theta\\ =&\left[\frac{3\pi}{64}- \frac{3\tan\frac{\theta}{2}}{160} - \tan^2\frac{\theta}{2}(\frac{3\pi}{32} - \frac{1}{256}) + \frac{3\pi\tan^4\frac{\theta}{2}}{64} +\frac{3\tan^5\frac{\theta}{2}}{160} - \frac{1}{768}\right]/(\tan^2\frac{\theta}{2} - 1)^3 -\\ &\frac{3\tanh^{-1}(\tan\frac{\theta}{2})}{160}.
\end{align*}
Plugging in the exact integral range, we thus obtain 
\begin{align*}
    &\int_{\cot^{-1}(\frac{1}{2})}^{\cot^{-1}(\frac{1}{2\sqrt{2}})}4\cos\theta\sin\theta\int_{-\pi/4}^{\pi/4}h_1(\theta,\phi)d\phi d\theta\\
   &\approx-0.026.
\end{align*}
Thus the negative part of $F$ is bounded by -0.026.
\subsubsection{\textbf{Positive part of F}}
 The positive part of $F$ is
\begin{equation}\label{PF}
       \int_{\cot^{-1}(\frac{1}{2})}^{\cot^{-1}(\frac{1}{2\sqrt{2}})}4\sin\theta\cos\theta\left(\int^{\pi/4}_{\cos^{-1}(2\cot\theta)} g_3(\theta,\phi)d\phi+\int_{-     \pi/4}^{-\cos^{-1}(2\cot\theta)} g_3(\theta,\phi)d\phi\right)d\theta
\end{equation}
\begin{align}
       \leq&\int_{\cot^{-1}(\frac{1}{2})}^{\cot^{-1}(\frac{1}{2\sqrt{2}})}4\sin\theta\cos\theta(\int^{\pi/4}_{\cos^{-1}(2\cot\theta)} \int_{2/2\cos\theta}^{4/2\sin\theta\cos\phi}t^2\Tilde{F}^{\omega_1}_{13}(t)dtd\phi\nonumber\\
       &+\int_{-     \pi/4}^{-\cos^{-1}(2\cot\theta)} \int_{2/2\cos\theta}^{4/2\sin\theta\cos\phi}t^2\Tilde{F}^{\omega_1}_{13}(t)dtd\phi)d\theta\label{pf1}.
\end{align}

But we observe that 
\begin{equation}\label{pf2}
    \int_{\cot^{-1}(\frac{1}{2})}^{\cot^{-1}(\frac{1}{2\sqrt{2}})}4\sin\theta\cos\theta\int_{-\cos^{-1}(2\cot\theta)}^{\cos^{-1}(2\cot\theta)}\int_{2/2\cos\theta}^{4/2\sin\theta\cos\phi}t^2\Tilde{F}^{\omega_1}_{13}(t)dtd\phi d\theta\geq 0,
\end{equation}

because when $\theta\in [\cot^{-1}(\frac{1}{2}),\cot^{-1}(\frac{1}{2\sqrt{2}})]$, $\phi\in [-\cos^{-1}(2\cot\theta),\cos^{-1}(2\cot\theta)]$,  
$$2/2\cos\theta>4/2\sin\theta\cos\phi,$$
and $\Tilde{F}^{\omega_1}_{13}(t)\leq 0$ when $t\leq 2/2\cos\theta$.

Hence (\ref{PF}) is bounded by (\ref{pf1})+(\ref{pf2}) which is
$$\int_{\cot^{-1}(\frac{1}{2})}^{\cot^{-1}(\frac{1}{2\sqrt{2}})}4\sin\theta\cos\theta\int_{-\pi/4}^{\pi/4}\int_{2/2\cos\theta}^{4/2\cos\theta\sin\theta}t^2\Tilde{F}^{\omega_1}_{13}(t)dtd\phi d\theta.$$
It has closed form
\begin{align*}
    &\int 4\sin\theta\cos\theta \int_{-\pi/4}^{\pi/4}\int_{2/2\cos\theta}^{4/2\cos\theta\sin\theta}t^2\Tilde{F}^{\omega_1}_{13}(t)dtd\phi d\theta\\
    =&8\tan\frac{\theta}{2}\frac{\log(\sqrt{2} + 1)+\sqrt{2}}{3}- \frac{32\log(\tan\frac{\theta}{2})}{3}\\
    +& \left[\tan^2\frac{\theta}{2}\frac{\pi+16}{6} + 8(\tan\frac{\theta}{2} - \tan^3\frac{\theta}{2})\frac{\log(\sqrt{2} + 1)+\sqrt{2}}{3} - \frac{8}{3}\right])/(\tan^2\frac{\theta}{2} - \tan^4\frac{\theta}{2})\\
    +& \frac{8}{3}\tan^2\frac{\theta}{2}.\\
\end{align*}
Plugging in the exact integral range, we thus obtain 
\begin{align*}
        &\int_{\cot^{-1}(\frac{1}{2})}^{\cot^{-1}(\frac{1}{2\sqrt{2}})}4\sin\theta\cos\theta \int_{-\pi/4}^{\pi/4}\int_{2/2\cos\theta}^{4/2\cos\theta\sin\theta}t^2\Tilde{F}^{\omega_1}_{13}(t)dtd\phi d\theta\\
    \approx&0.0064.
\end{align*}
Thus the positive part of $F$ is bounded by 0.0064.
\subsection{part G in case 3}
Part $G$ is 
\begin{align*}
    &\int_{\cot^{-1}\frac{3}{4}}^{\cot^{-1}\frac{3}{4\sqrt{2}}}4\sin\theta\cos\theta\int_{-\cos(\frac{3}{4}\cot\theta)}^{\cos(\frac{3}{4}\cot\theta)}h_1(\theta,\phi)+h_2(\theta,\phi)+g_3(\theta,\phi)d\theta\\
    &\int_{\cot^{-1}\frac{3}{4}}^{\cot^{-1}\frac{3}{4\sqrt{2}}}4\sin\theta\cos\theta(\int^{\pi/4}_{\cos^{-1}(\frac{3}{4}\cot\theta)}h_1(\theta,\phi)+h_2(\theta,\phi)+h_3(\theta,\phi)+g_4(\theta,\phi)d\phi\\
    &+\int_{-\pi/4}^{-\cos^{-1}(\frac{3}{4}\cot\theta)}h_1(\theta,\phi)+h_2(\theta,\phi)+h_3(\theta,\phi)+g_4(\theta,\phi)d\phi)d\theta,
\end{align*}

where $h_1,h_2$ are negative and $h_3,g_3,g_4$ are positive.
\subsubsection{\textbf{Negative part of G}}
The negative part of G is 
\begin{align*}
    &\int_{\cot^{-1}(\frac{3}{4})}^{\cot^{-1}(\frac{3}{4\sqrt{2}})}4\cos\theta\sin\theta\int_{-\pi/4}^{\pi/4}h_1(\theta,\phi)+h_2(\theta,\phi)d\phi d\theta\\
    =&\left(\frac{7\pi}{32}- \frac{9\tan\frac{\theta}{2}}{64} - \tan^2\frac{\theta}{2}(\frac{7\pi}{16} - \frac{31}{640}) + \frac{7\pi\tan^4\frac{\theta}{2}}{32}+ \frac{9\tan^5\frac{\theta}{2}}{64}- \frac{31}{1920}\right)/(\tan^2\frac{\theta}{2} - 1)^3\\ - &\frac{9\tanh^{-1}(\tan\frac{\theta}{2})}{64}\\
    \approx&-0.0694.
\end{align*}
Thus the negative part of $G$ is bounded by -0.0694.
\subsubsection{\textbf{Positive part of G}} The positive part of G can be split into two terms.
\begin{align*}
    &\textit{(G1)}\int^{\cot^{-1}(\frac{3}{4\sqrt{2}})}_{\cot^{-1}(\frac{3}{4})}4\cos\theta\sin\theta\int^{\cos^{-1}(\frac{3}{4}\cot\theta)}_{-\cos^{-1}(\frac{3}{4}\cot\theta)}g_3(\theta,\phi)d\phi d\theta\\
    &+\textit{(G2)} \int^{\cot^{-1}(\frac{3}{4\sqrt{2}})}_{\cot^{-1}(\frac{3}{4})}4\cos\theta\sin\theta(\int^{\pi/4}_{\cos^{-1}(\frac{3}{4}\cot\theta)}h_3(\theta,\phi)+g_4(\theta,\phi)d\phi+\int_{-\pi/4}^{-\cos^{-1}(\frac{3}{4}\cot\theta)}h_3(\theta,\phi)+g_4(\theta,\phi)d\phi)d\theta.
\end{align*}

Now
\begin{align}
    \textit{(G1)}&\int^{\cot^{-1}(\frac{3}{4\sqrt{2}})}_{\cot^{-1}(\frac{3}{4})}4\cos\theta\sin\theta\int^{\cos^{-1}(\frac{3}{4}\cot\theta)}_{-\cos^{-1}(\frac{3}{4}\cot\theta)}g_3(\theta,\phi)d\phi\nonumber\\
    \leq& \int^{\cot^{-1}(\frac{3}{4\sqrt{2}})}_{\cot^{-1}(\frac{3}{4})}4\cos\theta\sin\theta\int^{\cos^{-1}(\frac{3}{4}\cot\theta)}_{-\cos^{-1}(\frac{3}{4}\cot\theta)}\int_{2/2\cos\theta}^{4/2\sin\theta\cos\phi}t^2\Tilde{F}^{\omega_1}_{13}(t)f^{\omega_2}(t)f^{\omega_3}(t)dtd\phi,\label{pg1}
\end{align}
and
\begin{align*}
    \textit{(G2)}&\int^{\cot^{-1}(\frac{3}{4\sqrt{2}})}_{\cot^{-1}(\frac{3}{4})}4\cos\theta\sin\theta\left(\int^{\pi/4}_{\cos^{-1}(\frac{3}{4}\cot\theta)}h_3(\theta,\phi)+g_4(\theta,\phi)d\phi+\int_{-\pi/4}^{-\cos^{-1}(\frac{3}{4}\cot\theta)}h_3(\theta,\phi)+g_4(\theta,\phi)d\phi\right)
\end{align*}
is bounded by
\begin{align}\label{pg2}
    \int^{\cot^{-1}(\frac{3}{4\sqrt{2}})}_{\cot^{-1}(\frac{3}{4})}4\cos\theta\sin\theta(\int^{\pi/4}_{\cos^{-1}(\frac{3}{4}\cot\theta)}\int_{2/2\cos\theta}^{4/2\sin\theta\cos\phi}t^2\Tilde{F}^{\omega_1}_{13}(t)f^{\omega_2}(t)f^{\omega_3}(t)dtd\phi\nonumber\\
    +\int_{-\pi/4}^{-\cos^{-1}(\frac{3}{4}\cot\theta)}\int_{2/2\cos\theta}^{4/2\sin\theta\cos\phi}t^2\Tilde{F}^{\omega_1}_{13}(t)f^{\omega_2}(t)f^{\omega_3}(t)dtd\phi)
\end{align}
Hence the positive part of G is bounded by (\ref{pg1})+(\ref{pg2}), which has closed form
\begin{align*}
    &\int4\cos\theta\sin\theta\int_{-\pi/4}^{\pi/4}\int_{2/2\cos\theta}^{4/2\sin\theta\cos\phi}t^2\Tilde{F}^{\omega_1}_{13}(t)f^{\omega_2}(t)f^{\omega_3}(t)dtd\phi d\theta\\
    =&[\tan^6\frac{\theta}{2}\frac{\pi+2}{6} - \tan^4\frac{\theta}{2}\frac{10\pi+29}{30}+ \tan^2\frac{\theta}{2}\frac{15\pi+89}{90}+ \tan\frac{\theta}{2}\frac{10\log(\sqrt{2} + 1)+2\sqrt{2}+4}{15}\\
    +& \tan^5\frac{\theta}{2}(2\log(\sqrt{2} + 1) + \frac{2\sqrt{2}+4}{5}) - \tan^3\frac{\theta}{2}(2\log(\sqrt{2}+ 1) +\frac{+4\sqrt{2}+9}{10})\\
    -& \tan^7\frac{\theta}{2}(\frac{2\log(\sqrt{2} + 1)}{3}+ \frac{2\sqrt{2}}{15} + \frac{1}{6}) - \frac{1}{3}]/\left[\tan^2\frac{\theta}{2} - 3\tan^4\frac{\theta}{2} + 3\tan^6\frac{\theta}{2} - \tan^8\frac{\theta}{2}\right]\\ -&\frac{4\log(\tan\frac{\theta}{2})}{3} - \frac{\tanh\left(\frac{6391}{9\tan\frac{\theta}{2} + 240} - \frac{80}{3}\right)}{10} + \tan\frac{\theta}{2}\left(\frac{2\log(\sqrt{2} + 1)}{3} + \frac{2\sqrt{2}+4}{15}\right) + \frac{\tan^2\frac{\theta}{2}}{3}.
\end{align*}
Plugging in the exact integral range, we thus obtain
\begin{align*}
    &\int_{\cot^{-1}(\frac{3}{4})}^{\cot^{-1}(\frac{3}{4\sqrt{2}})}4\cos\theta\sin\theta\int_{-\pi/4}^{\pi/4}\int_{2/2\cos\theta}^{4/2\sin\theta\cos\phi}t^2\Tilde{F}^{\omega_1}_{13}(t)f^{\omega_2}(t)f^{\omega_3}(t)dtd\phi d\theta\\
    \approx& 0.0139.
\end{align*}
Therefore the positive part of $G$ is bounded by 0.0139.
\subsection{Altogether}
The sum of all further estimate cases is 
\begin{align*}
    &\textit{(Negative part of case 2+4)}+\textit{(Positive part of case 2+4)}+\textit{(Negative part of F)}\\
    &+\textit{(Positive part of F)}+\textit{(Negative part of G)}+\textit{(Positive part of G)}\\
    \leq& -0.0607+0.08718+(-0.026)\\
    &+0.0064+(-0.0694)+0.0139
    <0.
\end{align*}

\begin{remark}
It is very likely to extend the ideas to all dimensions $n \geq 4$ and show the correspondent integral is negative. 
\end{remark}

\end{CJK}
\end{document}